\documentclass[10pt,reqno]{amsart}
\usepackage{xifthen}
\usepackage{subfig}
\usepackage{pkgfile}

\newcommand{\Log}{{\rm Log}}
\begin{document}

\title[Limiting spectral measure for sums of unitary \& orthogonal matrices]{Limiting spectral distribution of \\ sums of unitary and orthogonal matrices}

\author[A.\ Basak]{Anirban Basak$^*$}
\thanks{${}^*$Research supported by Melvin and Joan Lane endowed Stanford Graduate Fellowship Fund.}
 \author[A.\ Dembo]{Amir Dembo$^\dagger$}\thanks{${}^\dagger$Research partially supported by NSF grant DMS-1106627.}
 
 \address{$^\dagger$Department of Mathematics, Stanford University
\newline\indent Building 380, Sloan Hall, Stanford, California 94305}

\address{$^{*\dagger}$Department of Statistics, Stanford University
\newline\indent Sequoia Hall, 390 Serra Mall, Stanford, California 94305}

\date{\today}



\maketitle

\begin{abstract}
We show that the empirical spectral distribution 
for sum of $d$ independent Haar distributed 
$n$-dimensional unitary matrices, converge for $n \to \infty$ to the Brown measure of the 
free sum of $d$ Haar unitary operators. The same applies for independent Haar distributed 
$n$-dimensional orthogonal matrices. As a byproduct of our approach, we relax the 
requirement of uniformly bounded imaginary part of Stieltjes transform of $T_n$
that is made in \cite[Theorem 1]{gkz}.
\end{abstract}

\section{Introduction}\label{sec-intro}
The method of moments and the Stieltjes transform approach 
provide rather precise information on asymptotics of the 
Empirical Spectral Distribution (in short \abbr{ESD}), 
for many {\em Hermitian} random matrix models.
In contrast, both methods fail for {\em non-Hermitian} matrix models, 
and the only available general scheme for finding the limiting spectral 
distribution in such cases is the one proposed by Girko (in \cite{girko}). 
It is extremely challenging to rigorously justify this scheme, even 
for the matrix model consisting of i.i.d. entries (of zero mean and finite variance). 
Indeed, after rather long series of partial results 
(see historical references in \cite{bordenave}), the {\em circular law} conjecture, for the i.i.d. case, 
was only recently established by Tao and Vu \cite{tao_vu} in full generality. 
Barring this simple model, very few results are known in the non-Hermitian regime. 
For example, nothing is known about the spectral 
measure of {\em random oriented $d$-regular} graphs. 
In this context, it was recently 
conjectured in \cite{bordenave} that, for $d \ge 3$, 
the \abbr{ESD} for the adjacency matrix of a uniformly chosen 
random oriented $d$-regular graph
converges to a measure $\mu_d$ on the complex plane, whose density
with respect to Lebesgue measure $m(\cdot)$ on $\C$ is  
\beq\label{eq:dreg}
h_d(v) := \f{1}{\pi} \f{d^2(d-1)}{(d^2 - |v|^2)^2}\bI_{\{|v| \le \sqrt{d}\}} \,.
\eeq
This conjecture, due to the observation that $\mu_d$ is 
the {\em Brown measure} of the {\em free sum} of $d \ge 2$ 
{\em Haar unitary operators} (see 
\cite[Example 5.5]{haagerup_larsen}), motivated us to consider the 
related problem of sum of $d$ independent Haar distributed, 
unitary or orthogonal matrices, for which we prove 
such convergence of the \abbr{ESD} in 
Theorem \ref{thm:orthogonal}. To this end, using hereafter 
the notation $\langle \Log, \mu \rangle_a^b := \int_a^b \log |x| d\mu(x)$ 
for any $a<b$ and probability measure $\mu$ on $\R$ 
(for which such integral is well defined), with 
$\langle \Log, \mu \rangle := \int_{\R} \log |x| d\mu(x)$, we
first recall the definition of Brown measure for a bounded 
operator (see \cite[Page 333]{haagerup_larsen}, or \cite{biane,brown}). 
\begin{dfn}
Let $(\cA, \tau)$ be a non-commutative $W^*$-probability space, i.e. 
a von Neumann algebra $\cA$ with a normal faithful tracial 
state $\tau$ (see \cite[Defn. 5.2.26]{agz}). 
For $h$ a positive element in $\cA$, let $\mu_h$ denote the unique 
probability measure on $\R^+$ such that 
$\tau(h^n) = \int t^n d \mu_h(t)$ for all $n \in \Z^+$. 
The Brown measure $\mu_a$ associated with each bounded $a \in \cA$, 
is the Riesz measure corresponding to the 
$[-\infty,\infty)$-valued sub-harmonic function 
$v \mapsto  \langle \Log, \mu_{|a-v|} \rangle$ on $\C$.   
That is, $\mu_a$ is the unique Borel probability measure on $\C$ such that   
\beq\label{def:bmeas}
d \mu_a (v) = \f{1}{2 \pi} \Delta_v 
\langle \Log, \mu_{|a-v|} \rangle 
\, dm(v)
, 
\eeq
where $\Delta_v$ denotes the two-dimensional Laplacian operator 
(with respect to $v \in \C$), and the identity (\ref{def:bmeas}) 
holds in distribution sense (i.e. when integrated against 
any test function $\psi \in C^\infty_c(\C)$).
\end{dfn}
\begin{thm}\label{thm:orthogonal}
For any $d\ge 1$, and $0 \le d' \le d$, as $n \ra \infty$ the 
\abbr{ESD} for sum of $d'$ independent, 
Haar distributed, $n$-dimensional unitary matrices $\{U_n^i\}$, and $(d-d')$ independent, Haar distributed, 
$n$-dimensional orthogonal matrices $\{O_n^i\}$, converges weakly, in probability, to the Brown measure $\mu_d$ of the
free sum of $d$ Haar unitary operators (whose density is 
given in (\ref{eq:dreg})). 
\end{thm}

\noindent
Recall that as $n \to \infty$, independent Haar distributed $n$-dimensional unitary (or orthogonal) 
matrices converge in $\star$-moments (see \cite{sniady} for a definition), to the 
collection $\{u_i\}_{i =1}^d$ of $\star$-{\em free} Haar unitary operators 
(see \cite[Theorem 5.4.10]{agz}). However, convergence of $\star$-moments, or even the stronger
convergence in distribution of traffics (of \cite{male}), do not necessarily imply convergence 
of the corresponding Brown measures\footnote{The Brown measure of a matrix 
is its \abbr{ESD} (see \cite[Proposition 1]{sniady})} 
(see \cite[\textsection 2.6]{sniady}). While \cite[Theorem 6]{sniady} shows that if the original matrices 
are perturbed by adding small Gaussian (of {\em unknown variance}), then the Brown measures do converge,
removing the Gaussian, or merely identifying the variance needed, are often hard tasks. 
For example, \cite[Prop. 7 and Cor. 8]{gwz} provide an example of ensemble where no
Gaussian matrix of polynomially vanishing variance can regularize the Brown measures (in this sense). 
Theorem \ref{thm:orthogonal} shows that sums of independent Haar distributed 
unitary/orthogonal matrices are {\em smooth} enough to have the 
convergence of \abbr{ESD}-s to the 
corresponding Brown measures {\em without adding any} Gaussian.

\noindent
Guionnet, Krishnapur and Zeitouni show in \cite{gkz} that the 
limiting \abbr{ESD} of $U_n T_n$ for non-negative definite, 
diagonal $T_n$ of limiting spectral measure $\Theta$, that is independent 
of the Haar distributed unitary (or orthogonal) matrix $U_n$, exists, is 
supported on a {\em single ring} and given by the Brown measure of the 
corresponding bounded (see \cite[Eqn. (1)]{gkz}), limiting operator. 
Their results, as well as our work, follow Girko's method, which 
we now describe, in brief.

\noindent
From Green's formula, for any polynomial 
$P(v)= \prod_{i=1}^n (v - \lambda_i)$ and test function 
$\psi \in C_c^2(\C)$, we have that 
\beq
\sum_{j=1}^n \psi(\lambda_j) = 
\f{1}{2 \pi} \int_\C \Delta \psi(v) \log |P(v)| dm(v) \,. 
\notag
\eeq
Considering this identity for the characteristic polynomial 
$P(\cdot)$ of a matrix $S_n$ (whose \abbr{ESD} we 
denote hereafter by $L_{S_n}$), results with 
\begin{align*}
\int_\C \psi(v) dL_{S_n}(v) =  & 
\f{1}{2 \pi n} \int_\C \Delta \psi(v) \log | \det  (vI_n-S_n)| dm(v)\\
= & \f{1}{4 \pi n} \int_\C \Delta \psi(v) \log  
\det  [(vI_n-S_n)(vI_n-S_n)^*] dm(v).
\end{align*} 
Next, associate with any $n$-dimensional 
non-Hermitian matrix $S_n$ and every $v \in \C$ 
the $2n$-dimensional {\em Hermitian} matrix 
\beq\label{eq:herm-def}
H_n^v:= 
 \begin{bmatrix}
  0 & (S_n -v I_n) \\
  (S_n - vI_n)^* & 0   
 \end{bmatrix} \,.
\eeq  
It can be easily checked that the eigenvalues of $H_n^v$ are merely
$\pm 1$ times the singular values of $vI_n- S_n$. Therefore, with $\nu_n^v$ 
denoting the \abbr{ESD} of $H_n^v$, we have that
\beq
\f{1}{n} \log \det[(vI_n-S_n)(vI_n-S_n)^*] = \f{1}{n} \log |\det H_n^v| 
= 2 \langle \Log, \nu_n^v \rangle \, , \notag 
\eeq
out of which we deduce the key identity
\beq
\int_\C \psi(v)dL_{S_n}(v) = \f{1}{2 \pi} 
\int_\C \Delta \psi(v) \langle \Log, \nu_n^v \rangle dm(v) 
\label{eq:girko_key_identity} 
\eeq
(commonly known as Girko's formula). 
The utility of Eqn. (\ref{eq:girko_key_identity}) lies in the following 
general recipe for proving convergence of $L_{S_n}$ per  
given family of non-Hermitian random matrices $\{S_n\}$ 
(to which we referred already as Girko's method).
\vskip5 pt

\noindent
 {\bf Step 1}: Show that for (Lebesgue almost) 
every $v \in \C$, as $n \ra \infty$ the measures 
$\nu_n^v$ converge weakly, in probability, to 
some measure $\nu^v$. 

\vskip5pt

\noindent
{\bf Step 2}: Justify that 
$\langle \Log, \nu_n^v \rangle \to \langle \Log, \nu^v \rangle$ 
in probability (which is the main technical challenge of this approach).

\vskip5pt

\noindent
{\bf Step 3}: A uniform integrability argument allows one to convert 
the $v$-a.e. convergence of $\langle \Log, \nu_n^v \rangle$ to 
the corresponding convergence for a suitable collection 
$\cS \subseteq C_c^2(\C)$ of (smooth) test functions. 
Consequently, it then follows from (\ref{eq:girko_key_identity}) that
for each fixed, non-random $\psi \in \cS$, 
\beq
\int_\C \psi(v) dL_{S_n}(v) \ra \f{1}{2 \pi} \int_\C \Delta \psi(v) 
\langle \Log, \nu^v \rangle
dm(v)\, , 
\label{eq:step3}
\eeq
in probability.

\vskip5pt
\noindent
{\bf Step 4}: Upon checking that $f(v) := \langle \Log, \nu^v \rangle$
is smooth enough to justify the integration by parts, one has that 
for each fixed, non-random $\psi \in \cS$, 
\beq
\label{eq:step4}
\int_\C \psi(v) dL_{S_n}(v)  \ra \f{1}{2 \pi}\int_\C 
\psi(v) \Delta f(v)
dm(v) \,, 
\eeq
in probability. For $\cS$ large enough, this implies the 
convergence in probability of the \abbr{ESD}-s 
$L_{S_n}$ to a limit which has the density 
$\frac{1}{2\pi} \Delta f$ 
with respect to Lebesgue measure on $\C$.

\vskip5pt

\noindent
Employing this method in \cite{gkz} requires, for 
{\bf Step 2}, to establish
suitable asymptotics for singular values of $T_n + \rho U_n$. Indeed, 
the key to the proofs there is to show that 
uniform boundedness of the imaginary part of 
Stieltjes transform of $\bm{T}_n$ (of the form assumed in \cite[Eqn. (3)]{gkz}), 
is inherited by the corresponding transform of $\bm{T}_n + \rho \bm{U}_n$ 
(see (\ref{eq:Y_n_T_n}) for a definition of $\bm{U}_n$ and $\bm{T}_n$). 
In the context of Theorem \ref{thm:orthogonal} (for $d' \ge 1$), 
at the start $d=1$, the expected \abbr{ESD} for $|v I_n - U_n|$ 
has unbounded density (see Lemma \ref{lem:u_n_density}), so the 
imaginary parts of relevant Stieltjes transforms are {\em unbounded}. 
We circumvent this problem by localizing the techniques of \cite{gkz}, 
whereby we can follow the development of unbounded regions of the
resolvent via the map $\bm{T}_n \mapsto \bm{T}_n + \rho (\bm{U}_n+\bm{U}_n^*)$ 
(see Lemma \ref{lem:new}), so as to achieve the desired convergence of 
integral of the logarithm near zero, for Lebesgue almost every $z$. 
We note in passing that Rudelson and Vershynin showed in 
\cite{vershynin} that the condition of \cite[Eqn. (2)]{gkz} about 
minimal singular value can be dispensed off 
(see \cite[Cor. 1.4]{vershynin}), but the remaining uniform boundedness 
condition \cite[Eqn. (3)]{gkz} is quite rigid. For example, it excludes 
atoms in the limiting measure $\Theta$ (so does not allow even $T_n=I_n$, 
see \cite[Remark 2]{gkz}). As a by product of our work, we relax below 
this condition about Stieltjes transform of $T_n$ (compare (\ref{eq:modi_T_n}) 
with \cite[Eqn. (3)]{gkz}), thereby generalizing \cite[Theorem 1]{gkz}.
\begin{prop} \label{rmk:gkz_improved}
Suppose the \abbr{ESD} of $\R^+$-valued, diagonal matrices $\{T_n\}$ 
converge weakly, in probability, to some probability measure $\Theta$. 
Assume further that:
\begin{enumerate}
\item There exists finite constant $M$ so that
\beq
\lim_{n \ra \infty} \P(\|T_n \| >M)=0. 
\label{eq:assumption_1}
\eeq
\item  There exists a closed set $K \subseteq \R$ of zero Lebesgue 
measure such that for every $\vep >0$, some $\kappa_\vep > 0$, 
$M_\vep$ finite and all $n$ large enough,
\beq
\{z: \Im (z) > n^{-\kappa_\vep}, |\Im(G_{\bm{T}_n}(z))| > M_\vep \} 
\subset \{z : z \in \bigcup_{x \in K} B(x,\vep) \}\,,
\label{eq:modi_T_n}
\eeq
where $G_{\bm{T}_n}(z)$ is the Stieltjes transform of the 
symmetrized version of the \abbr{ESD} of $T_n$,
as defined in (\ref{eq:G_T_n}). 
\end{enumerate}
If $\Theta$ is not a (single) Dirac measure, then the following hold:
\begin{enumerate}[(a)]
\item The \abbr{ESD} of 
$A_n:=U_n T_n$ converges, in probability, to limiting 
probability measure $\mu_A$.

\item The measure $\mu_A$ possesses a radially-symmetric density 
$h_A(v) := \f{1}{2 \pi} \Delta_v \langle \Log, \nu^v \rangle$ 
with respect to Lebesgue measure on $\C$, where 
$\nu^v:= \t{\Theta} \boxplus \lambda_{|v|}$ 
is the free convolution 
(c.f. \cite[\textsection 5.3.3]{agz}), of
$\lambda_r = \f{1}{2} (\delta_r + \delta_{-r})$ and 
the symmetrized version $\t{\Theta}$ of $\Theta$.

\item The support of $\mu_A$ is single ring: There exists constants $0 \le a < b < \infty$ so that
\beq
\text{supp} \; \mu_A = \{ r e^{i \theta}: a \le r \le b\}. \notag
\eeq
Further, $a=0$ if and only if $\int x^{-2} d \Theta (x)= \infty$.

\item The same applies if $U_n$ is replaced by a Haar distributed orthogonal matrix $O_n$.
\end{enumerate} 
\end{prop}

\noindent
This extension accommodates $\Theta$ with atoms, unbounded density, or singular part, as long as 
(\ref{eq:modi_T_n}) holds (at the finite $n$-level). For example, 
Proposition \ref{rmk:gkz_improved} 
applies for $T_n$ diagonal having $[n p_i]$ entries equal $x_i$, for $p_i>0$, $i=1,2,\ldots,k \ge 2$, 
whereas the case of $T_n = \alpha I_n$ for some $\alpha >0$ is an immediate consequence 
of Theorem \ref{thm:orthogonal}.

\noindent
Our presentation of the proof of Theorem \ref{thm:orthogonal} starts with detailed argument 
for $d'=d$, namely, the sum of independent Haar distributed unitary matrices. That is, we
first prove the following proposition, deferring to Section \ref{sec:thm_ii} its extension 
to all $0 \le d' < d$.
\begin{prop}\label{thm:unitary}
For any $d \ge 1$, as $n \to \infty$ the \abbr{ESD} of sum of 
$d$ independent, Haar distributed, $n$-dimensional unitary matrices $\{U_n^i\}_{i=1}^d$, converges weakly, in probability, to the Brown measure $\mu_d$ 
of free sum of $d$ Haar unitary operators. 
\end{prop}

\noindent
To this end, for any $v \in \C$ and i.i.d. 
Haar distributed unitary matrices $\{U_n^i\}_{1 \le i \le d}$,
and orthogonal matrices $\{O_n^i\}_{1 \le i \le d}$, let
\begin{equation}\label{eq:def_V1}
\bm{U}^{1,v}_n :=
\begin{bmatrix}
0 & (U_n^1 - v I_n)\\
(U_n^1 -v I_n)^* & 0
\end{bmatrix}
,
\end{equation}
and define $\bm{O}^{1,v}_n$ analogously, with $O_n^1$ replacing 
$U_n^1$. Set $\bm{V}^{1,v}_n := \bm{U}^{1,v}_n$ if $d' \ge 1$
and $\bm{V}^{1,v}_n := \bm{O}^{1,v}_n$ if $d'=0$, then let
\begin{equation}\label{eq:def_Vd}
\bm{V}^{k,v}_n := \bm{V}^{k-1,v}_n + \bm{U}_n^k + (\bm{U}_n^k)^* 
:= \bm{V}^{k-1,v}_n +
\begin{bmatrix}
0 & U_n^k\\
0 & 0
\end{bmatrix}
+ \begin{bmatrix}
0 & 0\\
(U_n^k)^* & 0
\end{bmatrix}
, \text{ for } k=2,\ldots,d'\,,
\end{equation}
and replacing $\bm{U}_n^k$ by $\bm{O}_n^k$, continue similarly 
for $k=d'+1,\ldots,d$. Next, let $G_n^{d,v}$ denote the 
expected Stieltjes transform of $\bm{V}_n^{d,v}$. That is,
\beq\label{eq:gndv}
G_n^{d,v}(z) := \E\Big[ \f{1}{2n} \Tr (z 
I_{2n} - \bm{V}_n^{d,v})^{-1} \Big],
\eeq
where the expectation is over {\em all} relevant unitary/orthogonal 
matrices $\{U_n^i,\, O_n^i, i=1,\ldots,d\}$. Part (ii) of the 
next lemma, about the relation between unbounded regions 
of $G_n^{d,v}(\cdot)$, and $G_n^{d-1,v}(\cdot)$ 
summarizes the key observation leading to Theorem \ref{thm:orthogonal} 
(with part (i) of this lemma similarly leading to our improvement 
over \cite{gkz}). To this end, for any $\rho>0$ and {\em arbitrary} 
$n$-dimensional matrix $T_n$ (possibly random), which is independent 
of the unitary Haar distributed $U_n$, let
\beq \label{eq:Y_n_T_n}
\bm{Y}_n := \bm{T}_n + \rho (\bm{U}_n +\bm{U}_n^*) := 
\begin{bmatrix}
0 & T_n\\
T_n^* & 0
\end{bmatrix}
\,+ 
\rho
\begin{bmatrix}
0 & U_n\\
0 & 0
\end{bmatrix}
+ \rho 
\begin{bmatrix}
0 & 0\\
U_n^* & 0
\end{bmatrix}
\eeq
and consider the following two functions of $z \in \C^+$,  
\begin{align}
G_{\bm{T}_n} (z) &:= \f{1}{2n} \Tr (z I_{2n} - \bm{T}_n)^{-1}, \label{eq:G_T_n}
\\
G_n(z) &:= \E\Big[ \f{1}{2n} \Tr (z I_{2n} - \bm{Y}_n)^{-1} \, | \, \bm{T}_n \Big]\,. 
\label{eq:G_n} 
\end{align}
\begin{lem}\label{lem:new}
(i) Fixing $R$ finite, suppose that $\|T_n\| \le R$ and 
the \abbr{ESD} of $\bm{T}_n$ converges to
some $\t{\Theta}$. Then, there exist $0 < \kappa_1 < \kappa$ small enough, 
and finite $M_\vep \uparrow \infty$ as $\vep \decto 0$, 
depending only on $R$ and $\t{\Theta}$, 
such that for all $n$ large enough and $\rho \in [R^{-1},R]$, 
\begin{align}\label{eq:basic-relation}
\Im(z) > n^{-\kappa_1} \; \& \; |\Im(G_{n}(z))| > 2 M_\vep 
\;\; \Longrightarrow \;\; \exists \psi_n(z) \in \C^+, \;
& \; \Im(\psi_n(z)) > n^{-\kappa} \; \& \; 
|\Im(G_{\bm{T}_n} (\psi_n(z))| > M_\vep \notag \\
& \; \& \;  z - \psi_n(z) \in B(-\rho,\vep) \cup B(\rho,\vep) \,.
\end{align}
The same applies when $U_n$ is replaced by Haar orthogonal matrix $O_n$
(possibly with different values of $0<\kappa_1<\kappa$ and 
$M_\vep \uparrow \infty$).

\noindent
(ii) For any $R$ finite, $d\ge 2$ and $d' \ge 0$, there exist 
$0 < \kappa_1 < \kappa$ small enough and finite $M_\vep \uparrow \infty$, 
such that (\ref{eq:basic-relation}) continues to hold 
for $\rho=1$, all $n$ large enough, any $|v| \le R$ 
and some $\psi_n(\cdot) := \psi_n^{d,v}(\cdot) \in \C^+$, 
even when $G_n$ and $G_{\bm{T}_n}$, 
are replaced by $G_n^{d,v}$ and $G_n^{d-1,v}$, respectively.
\end{lem}

\noindent
Section \ref{sec:lem_new} is devoted to the 
proof of Lemma \ref{lem:new}, building on which we prove 
Proposition \ref{thm:unitary} in Section \ref{sec:unitary}. 
The other key ingredients of this proof, namely 
Lemmas \ref{lem:sing_val_conv} and \ref{lem:log_integrate}, 
are established in Section \ref{sec:section_proofs_lemmas}.
Finally, short outlines of the proofs of 
Theorem \ref{thm:orthogonal} and of Proposition 
\ref{rmk:gkz_improved}, are provided in Sections \ref{sec:thm_ii} 
and \ref{sec:prop}, respectively.

\section{Proof of Lemma \ref{lem:new}}
\label{sec:lem_new}

This proof uses quite a few elements from the proofs in \cite{gkz}. 
Specifically, focusing on the case of unitary matrices, once a 
particular choice of $\rho \in [R^{-1},R]$ 
and $\bm{T}_n$ is made in part (i), all the steps appearing in 
\cite[pp. 1202-1203]{gkz} carry through, so all the equations 
obtained there continue to hold here (with a slight modification 
of bounds on error terms in the setting of part (ii), as 
explained in the sequel). Since this part follows \cite{gkz}, 
we omit the details. It is further easy to check that the 
same applies for the estimates obtained in 
\cite[Lemma 11, Lemma 12]{gkz}, which are thus also used 
in our proof (without detailed re-derivation).

\noindent
{\em Proof of  (i)}: We fix throughout this proof 
a fixed realization of the matrix $\bm{T}_n$, 
so expectations are taken only over the randomness 
in the unitary matrix $U_n$. Having done so,
first note that from \cite[Eqn. (37)-(38)]{gkz} we get
\beq \label{eq:sd_n_eqn_T}
G_n(z) = G_{\bm{T}_n}(\psi_n(z)) - \wt{O}(n, z, \psi_n(z)) \,,
\eeq
for 
\beq
\psi_n(z):= z - \f{\rho^2 G_n(z)}{1+ 2 \rho G_{U}^n(z)} \;  ,
\label{eq:def_psi_n_T}
\eeq
and 
\beq
G_{U}^n (z) := \E\Big[ 
\f{1}{2n} \Tr \big\{\bm{U}_n(z I_{2n} - \bm{Y}_n)^{-1}\big\} \,|\, \bm{T}_n \Big]
\,, 
\notag
\eeq
where for all $z_1, z_2 \in \C^+$
\beq
\wt{O}(n,z_1,z_2) = \f{2 O(n,z_1,z_2)}{1+ 2 \rho G_{U}^n(z_1)}\, , 
\label{eq:o_til}
\eeq
with $O(n,z_1,z_2)$ as defined in \cite[pp. 1202]{gkz}. 
Thus, (\ref{eq:sd_n_eqn_T}) and (\ref{eq:def_psi_n_T}) provide a relation 
between $G_n$ and $G_{\bm{T}_n}$ which is very useful for our proof. 
Indeed, from \cite[Lemma 12]{gkz} we have that there exists a 
constant $C_1:=C_1(R)$ finite such that, for all large $n$, if 
$\Im (z) > C_1 n^{-1/4}$ then 
\beq
\Im (\psi_n(z)) \ge \Im(z)/2. 
\label{eq:lem12}
\eeq
Additionally, from \cite[Eqn. (34)]{gkz} we have that
\beq
\rho (G_n(z))^2 = 2 G_{U}^n(z) (1 + 2 \rho G_{U}^n(z)) - O_1(n,z)\,, 
\label{eq:relation_T1}
\eeq
where 
$O_1(\cdot, \cdot)$ is as defined in \cite[pp. 1203]{gkz}. To this end, 
denoting 
\beq \label{eq:define_R_T}
F(G_n(z)):= \f{\rho^2 G_n(z)}{1+ 2 \rho G_{U}^n(z)}\;,  
\eeq
and using (\ref{eq:relation_T1}), we obtain after some algebra 
the identity 
\beq
G_n(z) \Big[\rho^2 - F^2(G_n(z)) \Big] = F(G_n(z)) \Big[ 1 + \f{\rho O_1(n,z)}{1+2 \rho G_{U}^n(z)} \Big] \,. 
\label{eq:r^2_1_T}
\eeq
Since
\beq
1 + 2 \rho G_{U}^{n}(z) = \f{1}{2} \Big(1 + \sqrt{1 + 4 \rho^2 G_n(z)^2 + 4 \rho O_1(n,z)}\Big), \label{eq:relation_T2}
\eeq
where the branch of the square root is uniquely determined 
by analyticity and the known behavior of $G_{U}^{n}(z)$ and $G_n(z)$ as
$|z| \ra \infty$ (see \cite[Eqn. (35)]{gkz}), we further have that
\begin{align}
F(G_n(z))  & =  \f{2 \rho^2 G_n(z)}{1+ \sqrt{1 + 4 (\rho G_n(z))^2 + 4 \rho O_1(n,z)}} \notag\\
 & = \f{1}{2} \Big[ \f{\rho^2 G_n(z) \sqrt{1+ 4 (\rho G_n(z))^2+4 \rho O_1(n,z)}}{(\rho G_n(z))^2 
+  \rho O_1(n,z)} - \f{\rho^2 G_n(z)}{(\rho G_n(z))^2 +  \rho O_1(n,z)} \Big] \,. \label{eq:bounding_r_T}
\end{align}
The key to our proof is the observation that if 
$|\Im(G_n(z))| \to \infty$ and 
$O_1(n,z)$ remains small, then from
(\ref{eq:bounding_r_T}), and (\ref{eq:def_psi_n_T}) necessarily 
$F(G_n(z)) = z - \psi_n(z) \to \pm \rho$. So, if 
$\widetilde{O}(n,z,\psi_n(z))$ remains bounded then 
by (\ref{eq:sd_n_eqn_T}) also $|\Im(G_{\bm{T}_n}(\psi_n(z)))| \to \infty$,
yielding the required result.

\noindent
To implement this, fix $M=M_\vep \ge 10$ 
such that $6 M^{-1}_\vep \le \vep^2$ and recall
that by \cite[Lemma 11]{gkz} there exists finite constant $C_2:=C_2(R)$ 
such that, for all large $n$, if $\Im(z) >C_1 n^{-1/4}$ then 
\beq 
\label{eq:lem11}
 |1+ 2 \rho G_{U}^n(z)| > C_2 \rho [ \Im(z)^3 \wedge 1].
\eeq
Furthermore, we have (see \cite[pp. 1203]{gkz}),
\beq\label{eq:onz1z2}
 |O(n,z_1,z_2)| \le \f{C \rho^2}{n^2 |\Im (z_2)| \Im (z_1)^2 (\Im(z_1) \wedge 1)}\,. 
\eeq
Therefore, enlarging $C_1$ as needed, by (\ref{eq:o_til}), (\ref{eq:lem12}), 
and 
(\ref{eq:lem11}) we obtain that, for all large $n$,
\beq
|\wt{O}(n,z,\psi_n(z))| \le \f{C \rho} 
{n^2 |\Im(\psi_n(z))| \Im(z)^2 (\Im(z)^4 \wedge 1)} \le M_\vep \notag
\eeq
whenever $\Im(z) > C_1 n^{-1/4}$. This, together with (\ref{eq:sd_n_eqn_T}), 
shows that if $|\Im(G_n(z))|>2 M_\vep$, then 
$|\Im (G_{\bm{T}_n}(\psi_n(z)))| >M_\vep$. 
Now, fixing $0 < \kappa_1 < \kappa <1/4$ 
we get from (\ref{eq:lem12}) that $\Im (\psi_n (z)) > n^{-\kappa}$. 
It thus remains to show only that 
$F(G_n(z)) \in B(-\rho,\vep) \cup B(\rho,\vep)$.
To this end, note that 
\beq
|O_1(n,z)| \le \f{C \rho^2}{n^2 \Im(z)^2 (\Im(z) \wedge 1)} \label{eq:o_1}
\eeq
(c.f. \cite[pp. 1203]{gkz}).
Therefore, $O_1(n,z)=o(n^{-1})$ whenever $\Im (z) > C_1 n^{-1/4}$, and so the rightmost term in (\ref{eq:bounding_r_T}) is 
bounded by $M_\vep^{-1}$ whenever $|\Im(G_n(z))| > 2M_\vep$. 
Further, when $\Im(z) > C_1 n^{-1/4}$, $|\Im(G_n(z))| > 2M_\vep$
and $n$ is large enough so $|O_1(n,z)| \le 1$, we have that
for any choice of the branch of the square root,
\begin{eqnarray*}
\Bigg| \f{ \rho G_n(z) \sqrt{1+ 4 (\rho G_n(z))^2+4 \rho O_1(n,z)}}{(\rho G_n(z))^2 +  \rho O_1(n,z)}\Bigg| \le \f{\sqrt{1+4 |\rho G_n(z)|^2+4 |\rho O_1(n,z)|}}{|\rho G_n(z)|-1} \le 4 \,,
\end{eqnarray*}
resulting with $|F(G_n(z))| \le 3 \rho$. Therefore, using (\ref{eq:lem11}) and (\ref{eq:o_1}), we get 
from (\ref{eq:r^2_1_T}) that if $\Im(z) > C_1 n^{-1/4}$ and $|\Im(G_n(z))| > 2 M_\vep$, then
\beq
\Big|F^2(G_n(z))-\rho^2 \Big| \le 6 |G_n(z)|^{-1} \le 6 M_\vep^{-1} \le \vep^2 \,. \notag
\eeq
In conclusion, 
$z - \psi_n(z)=F(G_n(z)) \in B(\rho,\vep) \cup B(-\rho,\vep)$, as stated. 
Further, upon modifying the values of $\kappa_1<\kappa$ and $M_\vep$, 
this holds also when replacing $U_n$ by a Haar distributed 
orthogonal matrix $O_n$. Indeed, the same analysis applies 
except for adding to $O(n,z_1,z_2)$ of \cite[pp. 1202]{gkz} a 
term which is uniformly bounded by $n^{-1} |\Im(z_2)|^{-1} 
(\Im(z_1) \wedge 1)^{-2}$ (see \cite[proof of Theorem 18]{gkz}),
and using in this case \cite[Cor. 4.4.28]{agz} to control the 
variance of Lipschitz functions of $O_n$ (instead of $U_n$).

\noindent
{\em Proof of (ii)}: Consider first the case of $d'=d$. Then,
setting $\rho=1$, $\bm{T}_n=\bm{V}_n^{d-1,v}$, 
and $\bm{Y}_n=\bm{V}_n^{d,v}$, one may check that following the 
derivation of \cite[Eqn. (37)-(38)]{gkz}, now with 
{\em all} expectations taken also {\em over} $\bf{T}_n$, we get
that 
\beq \label{eq:sd_n_eqn_U}
G_n^{d,v}(z) = G_{n}^{d-1,v}(\psi_n^{d,v}(z)) - \wt{O}(n, z, \psi_n^{d,v}(z)) 
\,,
\eeq
for some $K < \infty$ and all $\{z \in \C^+ : \Im(z) \ge K\}$, where 
\beq
\psi_n^{d,v} (z) := z - \f{ G_n^{d,v}(z)}{1+ 
2 G_{U_n}^{d,v}(z)} \;, 
\label{eq:def_psi_n_U}
\eeq
\beq
G_{U_n}^{d,v} (z) := \E\Big[ \f{1}{2n} 
\Tr \big\{\bm{U}_n^d(z I_{2n} - \bm{V}_n^{d,v})^{-1}\big\} \Big] \,,
\notag
\eeq
and for any $z_1,z_2 \in \C^+$,
\beq
\wt{O}(n,z_1,z_2) := \f{2 O(n,z_1,z_2)}{1+ 2G_{U_n}^{d,v}(z_1)}\,. 
\notag
\eeq
Next, note that for some $C < \infty$ and  
any $\C$-valued function $f_d(U_n^1,\ldots,U_n^d)$ of
i.i.d. Haar distributed $\{U_n^i\}$ 
\beq\label{eq:concen_U}
\E [ (f_d - \E[f_d])^2 ]  \le d C \| f_d \|^2_L \,,
\eeq
where $\|f_d\|_L$ 
denotes the relevant coordinate-wise Lipschitz norm, i.e.
$$
\| f_d \|_L := \max_{j=1}^d \sup_{U_n^1,\ldots,U_n^d, \wt{U}_n \ne U_n^j} 
\; \f{| f_d(U_n^1,\ldots,U_n^d) - 
f_d(U_n^1,\ldots,U_n^{j-1},\wt{U}_n,U_n^{j+1},\ldots) 
|}{\| U_n^j- \wt{U}_n \|_2} \;. 
$$
Indeed, we bound the variance of $f_d$ by the (sum of $d$) 
second moments of martingale differences 
$D_j f_d := {\bf E} [ f_d | U_n^1,\ldots,U_n^j]
- {\bf E}[f_d|U_n^1,\ldots, U_n^{j-1}]$. By the 
independence of 
$\{U_n^i\}$ and definition of $\|f_d\|_L$, 
conditional upon $(U_n^1,\ldots,U_n^{j-1})$, the $\C$-valued 
function $U_n^j \mapsto D_j f_d$ is Lipschitz of norm 
at most $\| f_d \|_L$ in the sense of \cite[Ineq. (4.4.31)]{agz}.
It then easily follows from the concentration inequalities of 
\cite[Cor. 4.4.28]{agz}, that the second moment of this 
function is at most $C \|f_d\|_L^2$  
(uniformly with respect to $(U_n^1,\ldots,U_n^{j-1})$).

\noindent
In the derivation of \cite[Lemma 10]{gkz}, the corresponding 
error term $O(n,z_1,z_2)$ is bounded by a sum of finitely 
many variances of Lipschitz functions of the form 
$\f{1}{2n} \Tr \{ H(U_n^d) \}$, each of which has
Lipschitz norm of order $n^{-1/2}$, hence controlled by 
applying the concentration inequality (\ref{eq:concen_U}).
We have here the same type of bound on $O(n,z_1,z_2)$, except 
that each variance in question is now with respect 
to some function $\f{1}{2n} \Tr \{ H(U_n^1,\ldots,U_n^d) \}$ 
having coordinate-wise Lipschitz norm of order $n^{-1/2}$
(and with respect to the joint law of the 
i.i.d. Haar distributed unitary matrices). 
Collecting all such terms, 
we get here instead of (\ref{eq:onz1z2}), the slightly worse bound
\beq
|O(n,z_1,z_2)| = O \bigg( \f{1}{n |\Im(z_2)| 
\Im(z_1)^2 (\Im(z_1) \wedge 1)^2 (\Im(z_2) \wedge 1)} \bigg)
\label{eq:O_bound_d}
\eeq
(with an extra factor $(\Im(z_2) \wedge 1)^{-1}$ due to 
the additional randomness in $(z_2 I_{2n} - {\bm T}_n)^{-1}$). 
Using the modified bound (\ref{eq:O_bound_d}), we proceed as in 
the proof of part (i) of the lemma, to first 
bound $\wt{O}(n,z,\psi_n^{d,v}(z))$, $O_1(n,z)$, and 
derive the inequalities replacing (\ref{eq:lem12}) 
and (\ref{eq:lem11}). Out of these bounds, we establish 
the stated relation (\ref{eq:basic-relation}) between $G_n^{d,v}$ 
and $G_n^{d-1,v}$ upon following the same route
as in our proof of part (i). Indeed, when doing so, 
the only effect of starting with (\ref{eq:O_bound_d})
instead of (\ref{eq:onz1z2})
is in somewhat decreasing the positive constants 
$\kappa_1, \, \kappa$, while increasing 
each of the finite constants $\{M_\vep, \vep >0\}$.

\noindent
Finally, with \cite[Cor. 4.4.28]{agz} applicable 
also over the orthogonal group, our proof of  
(\ref{eq:concen_U}) extends to any
$\C$-valued function 
$f_d(U_n^1,\ldots,U_n^{d'},O_n^{d'+1},\ldots,O_n^d)$ 
of independent Haar distributed unitary/orthogonal
matrices $\{U_n^i,O_n^i\}$. Hence, as 
in the context of part (i), the same argument applies 
for $0 \le d'<d$ (up to adding 
$n^{-1} |\Im(z_2)|^{-1} (\Im(z_1) \wedge 1)^{-2}$ 
to 
(\ref{eq:O_bound_d}), 
c.f. \cite[proof of Theorem 18]{gkz}).

\section{Proof of Proposition \ref{thm:unitary}}
\label{sec:unitary}

It suffices to prove Proposition \ref{thm:unitary} only for $d \ge 2$, 
since the easier case of $d=1$ has already been established in 
\cite[Corollary 2.8]{MM}. We proceed to do so via the four steps 
of Girko's method, as described in Section \ref{sec-intro}. 
The following two lemmas (whose 
proof is deferred to Section \ref{sec:section_proofs_lemmas}), take care
of {\bf Step 1} and {\bf Step 2} of Girko's method, respectively.
\begin{lem}\label{lem:sing_val_conv}
Let $\lambda_1 = \f{1}{2} (\delta_{-1}+ \delta_1)$ and 
$\Theta^{d,v} := \Theta^{d-1,v} \boxplus \lambda_1$ for  
all $d \ge 2$, starting at $\Theta^{1,v}$ which for $v \ne 0$ 
is the symmetrized version of the measure on $\R^+$ having 
the density $f_{|v|} (\cdot)$ of (\ref{eq:f_r_expression}), 
while $\Theta^{1,0} = \lambda_1$. Then, for each $v \in \C$ 
and $d \in \N$, the \abbr{ESD}-s $L_{\bm{V}_n^{d,v}}$ of the  
matrices $\bm{V}_n^{d,v}$ (see (\ref{eq:def_Vd})), converge 
weakly as $n \to \infty$, in probability, to $\Theta^{d,v}$.
\end{lem}

\begin{lem} \label{lem:log_integrate}
For any $d \ge 2$ and Lebesgue almost every $v \in \C$,
\beq
\langle \Log, L_{\bm{V}_n^{d,v}} \rangle 
\ra \langle \Log ,\Theta^{d,v} \rangle, \label{eq:log_integrates}
\eeq
in probability. Furthermore, there exist closed $\Lambda_d \subset \C$ 
of zero Lebesgue measure, such that
\beq
\int_{\C} \phi(v) \langle \Log,   
L_{\bm{V}_n^{d,v}} \rangle  dm(v) \ra 
\int_{\C} \phi(v) \langle \Log,   \Theta^{d,v}\rangle dm (v), 
\label{eq:function_log_integrates}
\eeq
in probability, for each fixed, non-random $\phi \in C_c^\infty(\C)$ 
whose support is disjoint of $\Lambda_d$. That is, the support of 
$\phi$ is contained for some $\gamma >0$, in the bounded, open set 
\beq
\Gamma^d_{\gamma} := \big\{ v \in \C : \gamma < |v| < \gamma^{-1}, \; 
\inf_{u \in \Lambda_d} \{\, | v - u | \} > \gamma \big\} \,.
\label{eq:gamma_d}
\eeq
\end{lem}

\noindent
We claim that the convergence result of 
(\ref{eq:function_log_integrates}) provides us already 
with the conclusion (\ref{eq:step3}) 
of {\bf Step 3} in Girko's method, for test functions in
$$
\cS := \{ \psi \in C_c^\infty(\C),  \hbox{ supported within }
\Gamma^d_{\gamma} \hbox{ for some } \gamma > 0 
\} \,.
$$
Indeed, fixing $d \ge 2$, the Hermitian matrices
$\bm{V}_n^{d,v}$ of (\ref{eq:def_Vd}) are precisely 
those $H_n^v$ of the form (\ref{eq:herm-def}) 
that are associated with $S_n := \sum_{i=1}^d U_n^i$ in 
Girko's formula (\ref{eq:girko_key_identity}). Thus. 
combining the latter identity for $\psi \in \cS$ 
with the convergence result of (\ref{eq:function_log_integrates}) 
for $\phi = \Delta \psi$, we get the following 
convergence in probability as $n \to \infty$,
\beq
\int_{\C} \psi(v) d L_{S_n} (v) 
= \f{1}{2\pi}
\int_{\C} \Delta \psi (v) \langle \Log , 
L_{\bm{V}_n^{d,v}} \rangle   dm(v) 
\ra \f{1}{2 \pi} \int_{\C} \Delta \psi(v) \langle \Log,\Theta^{d,v}
\rangle d m(v) \,.
\label{eq:final_step}
\eeq

Proceeding to identify the limiting measure as the Brown measure 
$\mu_d := \mu_{s_d}$ of the sum $s_d:=u_1+u_2+\cdots+u_d$ 
of $\star$-free Haar unitary operators $u_i$, recall 
\cite{NS} that each $(u_i,u_i^*)$ is ${\bf R}$-diagonal.
Hence, by \cite[Proposition 3.5]{haagerup_larsen} 
we have that $\Theta^{d,v}$ is the symmetrized version of 
the law of $|s_d - v|$, and so by definition (\ref{def:bmeas}) 
we have that for any $\psi \in C_c^\infty (\C)$, 
\beq
\f{1}{2 \pi} \int_{\C} \Delta \psi(v) 
\langle \Log,\Theta^{d,v}\rangle d m(v) 
= \int_{\C} \psi(v) \mu_{s_d} (dv) \,.
\label{eq:bm-limit}
\eeq  
In parallel with {\bf Step 4} of Girko's method, it thus suffices for 
completing the proof, to verify that the convergence in probability
\beq\label{eq:func-conv}
\int_{\C} \psi(v) dL_{S_n}(v) \ra \int_{\C} \psi(v) d \mu_{s_d} (v) \,,
\eeq
for each {\em fixed} $\psi \in \cS$, yields the weak 
convergence, in probability, of $L_{S_n}$ to $\mu_{s_d}$.
 
To this end, suppose first that (\ref{eq:func-conv}) holds 
almost surely for each fixed $\psi \in \cS$, and recall that
for any $\gamma >0$ and each open $G \subset \Gamma_{\gamma}^d$
there exist $\psi_k \in \cS$ such that $\psi_k \uparrow 1_{G}$.
Consequently, a.s.
\beq
\liminf_{n \ra \infty} L_{S_n} (G) \ge \sup_k \liminf_{n \to \infty} 
\int_{\C} \psi_k(v) d L_{S_n} (v) = 
\sup_k \int_{\C} \psi_k (v) d\mu_{s_d}(v) = \mu_{s_d} (G) \,. \notag
\eeq
Further, from \cite[Example 5.5]{haagerup_larsen} we know that 
$\mu_{s_d}$ has, for $d \ge 2$, a bounded density 
with respect to Lebesgue measure on $\C$ (given by
$h_d(\cdot)$ of (\ref{eq:dreg})). In particular, since
$m(\Lambda_d)=0$, it follows that
$\mu_{s_d}(\Lambda_d) = 0$ and hence 
$\mu_{s_d}(\Gamma_\gamma^d) \to 1$ when $\gamma \to 0$. 
Given this, fixing some $\gamma_\ell \downarrow 0$ and 
open $G \subset \C$,  we deduce that a.s.
\beq
\liminf_{n \ra \infty} L_{S_n} (G) 
\ge \lim_{\ell \to \infty} \liminf_{n \to \infty} 
L_{S_n} (G \cap \Gamma_{\gamma_\ell}^d) 
\ge \lim_{\ell \ra \infty} \mu_{s_d} (G \cap \Gamma_{\gamma_\ell}^d)    
 = \mu_{s_d} (G)  \,.  
\label{eq:a.s._final}
\eeq
This applies for any countable collection $\{G_i\}$ 
of open subsets of $\C$, with the reversed inequality holding for 
any countable collection of closed subsets of
$\C$. In particular, fixing any countable convergence determining class 
$\{f_j\} \subset C_b(\C)$ and countable dense $\widehat{\Q} \subset \R$
such that $\mu_{s_d}(f_j^{-1}(\{q\}))= 0$ for all $j$ and $q \in \widehat{\Q}$, 
yield the countable collection $\cG$ of $\mu_{s_d}$-continuity sets
(consisting of interiors and complement of closures 
of $f_j^{-1}([q,q'))$, $q,q' \in \widehat{\Q}$), 
for which $L_{S_n}(\cdot)$ converges to 
$\mu_{s_d}(\cdot)$. The stated a.s. weak convergence of $L_{S_n}$ to 
$\mu_{s_d}$ then follows as in the usual proof of Portmanteau's 
theorem, under our 
assumption that (\ref{eq:func-conv}) holds a.s. 

\noindent
This proof extends to the case at hand, where (\ref{eq:func-conv}) 
holds in probability, since convergence in probability implies 
that for every subsequence, there exists a further subsequence along which 
a.s. convergence holds, and the whole argument uses only 
countably many functions $\psi_{k,\ell,i} \in \cS$. Specifically, 
by a Cantor diagonal argument, for any given subsequence 
$n_j$, we can extract a further 
subsequence $j(l)$, such that 
(\ref{eq:a.s._final}) holds a.s. for $L_{S_{n_{j(l)}}}$ and all
$G$ in the countable collection $\cG$ of $\mu_{s_d}$-continuity sets. 
Therefore, a.s. $L_{S_{n_{j(l)}}}$ 
converges weakly to $\mu_{s_d}$ 
and by the arbitrariness of $\{n_j\}$ we have that,
in probability, $L_{S_n}$ converges to $\mu_{s_d}$ weakly. 


\section{Proofs of Lemma \ref{lem:sing_val_conv} and Lemma \ref{lem:log_integrate}}
\label{sec:section_proofs_lemmas}

We start with a preliminary result, 
needed for proving Lemma \ref{lem:sing_val_conv}.
\begin{lem}\label{lem:u_n_density}
For Haar distributed $U_n$ and any $r>0$,
the expected \abbr{ESD} of $|U_n -r I_n|$
has the density 
\beq \label{eq:f_r_expression}
f_r(x)= \f{2}{\pi}\f{x}{\sqrt{(x^2- (r-1)^2)((r+1)^2- x^2)}} \,,  \ \ \ |r-1| \le x \le r+1
\eeq
with respect to Lebesgue's measure on $\R^+$ (while for $r=0$, 
this \abbr{ESD} consists of a single atom at $x=1$).
\end{lem}

\noindent
{\em Proof}: It clearly suffices 
to show that the expected \abbr{ESD} of  $(U_n -r I_n)(U_n -r I_n)^*$ has for $r>0$ 
the density
\beq\label{eq:grdef}
g_r(x)=\f{1}{\pi}\f{1}{\sqrt{(x- (r-1)^2)((r+1)^2 - x)}}, \ \ \ (r-1)^2 \le x \le (r+1)^2 \,.
\eeq
To this end note that by the invariance of the Haar unitary measure
under multiplication by $e^{i \theta}$, we have that
\beq\label{eq:unit-idn}
\E[\f{1}{n} \Tr \{U_n^k\}]= \E[\f{1}{n} \Tr \{(U_n^*)^k\}] = 0 \,,
\eeq 
for all positive integers $k$ and $n$. Thus,
\beq
\E\Big[ \f{1}{n} \Tr\big\{(U_n+U_n^*)^k\big\}\Big]= {k \choose k/2} 
\text{ for } k \text{ even and } 0 \text{ otherwise}.
\notag
\eeq
Therefore, by the moment method, the expected \abbr{ESD} of 
$U_n + U_n^*$ (denoted $\bar{L}_{U_n+U_n^*}$), satisfies
\beq
\bar{L}_{U_n+U_n^*} \stackrel{d}{=} 2 \cos \theta = e^{i \theta} 
+ e^{-i \theta}, \text{ where } 
\theta \sim \dU (0, 2 \pi).
\notag
\eeq
Consequently, we get the formula (\ref{eq:grdef}) 
for the density $g_r(x)$ of the expected \abbr{ESD} of 
\beq
(U_n - r I_n)(U_n - r I_n)^* = (1+r^2)I_n - r (U_n+U_n^*), \notag
\eeq
by applying the change of variable formula for $x=(1+r^2) - 2 r \cos \theta$
(and $\theta \sim \dU (0, 2 \pi)$).
\qed

\medskip
\noindent
{\em Proof of Lemma \ref{lem:sing_val_conv}}: 
Recall \cite[Theorem 2.4.4(c)]{agz} that for the claimed 
weak convergence of $L_{\bm{V}_n^{d,v}}$ to $\Theta^{d,v}$, 
in probability, it suffices to show that 
per fixed $z \in \C^+$, the corresponding Stieltjes 
transforms 
$$
f_n^{d,v} (z) := \f{1}{2n} \Tr \{ (z I_{2n} - \bm{V}_n^{d,v})^{-1} \} 
$$
converge in probability to the Stieltjes transform 
$G^{d,v}_\infty(z)$ of $\Theta^{d,v}$. To this end, 
note that each $f_n^{d,v} (z)$ is a point-wise 
Lipschitz function of $\{U_n^i\}$, whose expected 
value is $G_n^{d,v}(z)$ of (\ref{eq:gndv}),
and that $\|f_n\|_L \to 0$ as $n \to \infty$ (per fixed
values of $d,v,z$). It thus follows from (\ref{eq:concen_U}) 
that as $n \to \infty$,
\beq
\E[( f_n^{d,v} (z) - G_n^{d,v}(z))^2] \ra 0  
\notag
\eeq
and therefore, it suffices to prove that per fixed $d$, 
$v \in \C$ and $z \in \C^+$, as $n \to \infty$,  
\beq\label{eq:gnconv}
G_n^{d,v}(z) \to G_\infty^{d,v} (z) \,.
\eeq
Next observe that by invariance of the 
law of 
$U_n^1$
to multiplication by scalar $e^{i \theta}$, 
the expected \abbr{ESD} of 
$\bm{V}_n^{1,v}$ 
depends only on $r=|v|$,
with
$\Theta^{1,v} = \E [L_{\bm{V}_n^{1,v}}]$   
(see Lemma \ref{lem:u_n_density}). Hence, (\ref{eq:gnconv}) 
trivially holds for $d=1$ and we proceed to prove the latter
pointwise (in $z,v$), convergence by an induction on $d \ge 2$.
The key ingredient in the induction step is the (finite $n$) 
Schwinger-Dyson equation in our set-up, namely 
Eqn. (\ref{eq:sd_n_eqn_U})-(\ref{eq:def_psi_n_U}).  
Specifically, from (\ref{eq:sd_n_eqn_U})-(\ref{eq:def_psi_n_U}) 
and the induction hypothesis 
it follows that for some non-random $K<\infty$, 
any limit point, denoted $(G^{d,v}, G_U^{d,v})$, 
of the uniformly bounded, equi-continuous functions 
$(G_n^{d,v}, G_{U_n}^{d,v})$ on $\{z \in \C^+ : 
\Im(z) \ge K\}$, satisfies 
\beq
G^{d,v}(z)= G^{d-1,v}_\infty (\psi(z)), \text{ with } 
\psi(z):=z - \f{G^{d,v}(z)}{1+2 G_U^{d,v}(z)}\,.
\label{eq:sd_infinite}
\eeq
Moreover, from the equivalent version of (\ref{eq:relation_T1}) 
in our setting, we obtain that
\beq
4 G_U^{d,v}(z) = -1+ \sqrt{1+4 G^{d,v}(z)^2}\,, \notag
\eeq
for a suitable branch of the square root (uniquely determined 
by analyticity and decay to zero as $|z| \to \infty$ of 
$z \mapsto (G^{d,v}(z), G_{U}^{d,v}(z))$). 
Thus, $G(z) = G^{d,v}(z)$ satisfies the relation 
\beq
G(z) - G_\infty^{d-1,v} \Big( z - \f{2G(z)}{1+ \sqrt{1+ 4 G(z)^2}}
\Big) = 0 \,. 
\label{eq:psi_infinite}
\eeq
Since $\Theta^{d,v} = \Theta^{d-1,v} \boxplus \lambda_1$, it 
follows that (\ref{eq:psi_infinite}) holds also for 
$G(\cdot)=G_\infty^{d,v}(\cdot)$ (c.f. \cite[Remark 7]{gkz}). 
Further, $z \mapsto G_\infty^{d-1,v}(z)$ is analytic on $\C^+$ with 
derivative of $O(z^{-2})$ at infinity, hence by the implicit 
function theorem the identity (\ref{eq:psi_infinite}) uniquely 
determines the value of $G(z)$ for all $\Im (z)$ large enough. In particular, 
enlarging $K$ as needed, $G^{d,v}= G_\infty^{d,v}$ 
on $\{z \in \C^+: \Im (z) \ge K\}$, which by analyticity of both functions
extends to all of $\C^+$. With (\ref{eq:gnconv}) verified, this completes
the proof of the lemma. 
\qed

\medskip
\noindent
The proof of Lemma \ref{lem:log_integrate} requires 
the control of $\Im(G_n^{d,v}(z))$ as
established in Lemma \ref{lem:induction_step}. This is done
inductively in $d$, with Lemma \ref{lem:induction_start} 
providing the basis $d=1$ of the induction.
\begin{lem}\label{lem:induction_start}
For some $C$ finite, all $\vep \in (0,1)$ and $v \in \C$, 
\beq
\Big\{z \in \C^+  : \ |\Im G^{1,v}_n(z)| \ge C \vep^{-2}\Big\} 
\subseteq \Big \{E + i \eta : \eta \in (0,\vep^2), E \in 
\big(\pm(1 \pm |v|) - 2 \vep, \pm(1 \pm |v|)+ 2 \vep \big) \Big\}. 
\notag
\eeq
\end{lem}

\noindent
{\em Proof}: It is trivial to confirm our claim in case $v=0$ 
(as $G_n^{1,0}(z)=z/(z^2-1)$). Now, fixing $r=|v|>0$, let 
$\wt{f}_r(\cdot)$ denote the symmetrized version of the 
density $f_r(\cdot)$, and note that for any $\eta>0$,
\begin{eqnarray}
|\Im G^{1,v}_n (E+ i \eta)|  &= & \int_{|x-E| > \sqrt{\eta}} \f{\eta}{(x-E)^2 + \eta^2} \wt{f}_r (x) dx + \int_{|x-E| \le \sqrt{\eta}} \f{\eta}{(x-E)^2 + \eta^2} \wt{f}_r (x) dx \notag \\
& \le & 1+ \Big[\sup_{\{x: |x -E| \le \sqrt{\eta}\}} \wt{f}_r(x) \Big] \int_{|x-E| \le \sqrt{\eta}} \f{\eta}{(x-E)^2 + \eta^2} dx \notag \\
& \le & 1+ \pi \Big[\sup_{\{x: |x -E| \le \sqrt{\eta}\}} \wt{f}_r(x) \Big]. \label{eq:initi_bd_ImG}
\end{eqnarray}
With $\Gamma_\vep$ denoting the union of open intervals of radius $\vep$ 
around the four points $\pm 1 \pm r$, it follows 
from (\ref{eq:f_r_expression}) that for some $C_1$ finite 
and any $r,\vep>0$,
\beq
\sup_{x \notin \Gamma_\vep} \{\wt{f}_r (x)\} \le C_1 \vep^{-2}\,. \notag 
\eeq
Thus, from (\ref{eq:initi_bd_ImG}) it follows that
\beq
\sup_{\{E, \eta: (E - \sqrt{\eta}, E+ \sqrt{\eta}) \subset \Gamma_\vep^c\}} \; |\Im G^{1,v}_n(E+ i \eta)| \le C \vep^{-2}\,,
\notag
\eeq
for some $C$ finite, all $\vep \in (0,1)$ and $r>0$.
To complete the proof simply note that 
\beq
\{(E,\eta): E\in \Gamma_{2 \vep}^c, \eta \in (0,\vep^2)\} \subseteq
\{(E,\eta): (E-\sqrt{\eta}, E+\sqrt{\eta}) \subseteq \Gamma_\vep^c\}, \notag
\eeq
and
\beq
\sup_{E \in \R, \eta \ge \vep^2} |\Im G^{1,v}_n (E+ i \eta)| \le \vep^{-2}. 
\notag
\eqno{\qed}
\eeq

\medskip
\noindent
Since the density $\wt{f}_{|v|} (\cdot)$ is unbounded at $\pm 1 \pm |v|$, 
we can not improve Lemma \ref{lem:induction_start} to show 
that $\Im G_n^{1,v}(z)$ is uniformly bounded. 
The same applies for $d \ge 2$ so a result such as 
\cite[Lemma 13]{gkz} is not possible in our set-up. Instead, 
as we show next, inductively applying 
Lemma \ref{lem:new}(ii) allows us to control 
the region where $|\Im (G_n^{d,v}(z))|$ might 
blow up, in a manner which suffices for 
establishing Lemma \ref{lem:log_integrate} 
(and consequently Proposition \ref{thm:unitary}).
\begin{lem} \label{lem:induction_step}
For $r \ge 0$, $\gamma > 0$ and integer $d \ge 1$, let
$\Gamma_\gamma^{d,r} \subset \C$ denote the union of 
open balls of radius $\gamma$ centered at $\pm m \pm r$ 
for $m=0,1,2,\ldots,d$.
Fixing integer $d \ge 1$, $\gamma \in (0,1)$ and $R$ finite, 
there exist $M$ finite and $\kappa >0$ such that 
for all $n$ large enough and any $v \in B(0,R)$,
\beq \label{eq:unif_bound}
\sup\{ |\Im(G_n^{d,v}(z))| : \; \Im (z) > n^{-\kappa}, \; 
z \notin \Gamma_\gamma^{d,|v|} \} \le M \,. 
\eeq
\end{lem}

\noindent
{\em Proof}: For any $d \ge 1$, $v \in \C$, positive $\kappa$ and
finite $M$, set   
$$
\Gamma_n^{d,v}(M,\kappa) 
:= \{z: \Im(z) > n^{-\kappa}, |\Im(G_n^{d,v}(z))| >M \}
\,,
$$ 
so our thesis amounts to the existence of finite 
$M$ and $\kappa>0$, depending only on $R$, 
$d \ge 2$ and $\gamma \in (0,1)$, such that for all $n$ large enough,
\beq
\Gamma_n^{d,v}(M,\kappa) \subset \Gamma_\gamma^{d,|v|}, \qquad 
\qquad \forall v \in B(0,R) \,.
\label{eq:induction_hypothesis}
\eeq 
Indeed, for $d=1$ this is a direct consequence of Lemma 
\ref{lem:induction_start} (with $\gamma = 2\vep$, $M=C \vep^{-2}$), 
and we proceed to confirm (\ref{eq:induction_hypothesis}) 
by induction on $d \ge 2$. To carry out the inductive step from 
$d-1$ to $d$, fix $R$ finite and $\gamma \in (0,1)$, assuming that
(\ref{eq:induction_hypothesis}) applies at $d-1$ and $\gamma/2$, 
for some finite $M_\star$ and positive $\kappa_\star$ 
(both depending only on $d$, $R$ and $\gamma$). Then, let 
$\vep \in (0,\gamma/2)$ be small enough such that 
Lemma \ref{lem:new}(ii) applies for some 
$M_\vep \ge M_\star$ and $0 < \kappa_1 < \kappa \le \kappa_\star$. 
From Lemma \ref{lem:new}(ii) we know that for any $n$ large 
enough, $v \in B(0,R)$ and $z \in \Gamma_n^{d,v}(2 M_\vep,\kappa_1)$,
there exists $w := \psi_n^{d,v} (z)$ for which
$$
z - w \in B(-1,\vep) \cup B(1,\vep) 
\qquad \& \qquad  
w \in \Gamma_n^{d-1,v}(M_\vep,\kappa) \subseteq 
\Gamma_n^{d-1,v}(M_\star,\kappa_\star) \subset \Gamma_{\gamma/2}^{d-1,|v|} 
\,,
$$
where the last inclusion is due to our choice of $M_\star$ and $\kappa_\star$. 
With $\vep \le \gamma/2$, it is easy 
to check that $z - w \in B(-1,\vep) \cup B(1,\vep)$ and 
$w \in \Gamma_{\gamma/2}^{d-1,r}$ result with 
$z \in \Gamma_{\gamma}^{d,r}$. That is, we have established the validity of 
(\ref{eq:induction_hypothesis}) at $d$ and arbitrarily small $\gamma$, 
for $M=2 M_\vep$ finite and $\kappa_1$ positive, 
both depending only on $R$, $d$ and $\gamma$.
\qed

\vskip 5 pt
\noindent
{\em Proof of Lemma \ref{lem:log_integrate}}: 
Recall \cite[Theorem 1.1]{vershynin} the existence of universal constants 
$0< c_1$ and $c_2 < \infty$, such that for any non-random matrix $D_n$ 
and Haar distributed unitary matrix $U_n$, the smallest singular value 
$s_{\min}$ of $U_n+D_n$ satisfies,
\beq
\P(s_{\min}(U_n+D_n) \le t) \le t^{c_1} n^{c_2}. \label{eq:rv_unitary}
\eeq
The singular values of $\bm{V}_n^{d,v}$ are clearly the same as those 
of $S_n - v I_n = U_n^1 +D_n$ for $D_n=\sum_{i=2}^{d} U_n^i - v I_n$,
which is independent of the Haar unitary $U_n^1$. Thus, 
applying (\ref{eq:rv_unitary}) conditionally on $D_n$, we get that  
\beq
\P(s_{\min}(\bm{V}_n^{d,v}) \le t) \le t^{c_1} n^{c_2} \,, \label{eq:smin_U}
\eeq
for every $v \in \C$, $t>0$ and $n$. It then follows that 
for any $\delta>0$ and $\alpha < c_1$, 
\beq
\E \Big[ (s_{\min}(\bm{V}_n^{d,v}))^{-\alpha}  
\bI_{\big\{s_{\min}(\bm{V}_n^{d,v})\le n^{-\delta}\big\}}\Big]
\le \f{c_1}{c_1-\alpha} n^{c_2-\delta (c_1-\alpha)} \,.
\label{eq:smin_U2}
\eeq
Setting hereafter $\alpha=c_1/2$ positive and $\delta=4 c_2/c_1$ finite, 
the right side of (\ref{eq:smin_U2}) decays to zero as $n \to \infty$. 
Further, for any $n$, $d$ and $v$, 
\beq
\E \Big[
\langle |\Log|, L_{\bm{V}_n^{d,v}} \rangle_0^{n^{-\delta}} 
\Big] 
\le \E \Big[   |\log s_{\min}(\bm{V}_n^{d,v})| 
\bI_{\big\{s_{\min}(\bm{V}_n^{d,v})\le n^{-\delta}\big\}}\Big] \;.
\label{eq:smin_U1}
\eeq
Hence, with $|x|^\alpha \log |x| \to 0$ as $x \to 0$, upon 
combining (\ref{eq:smin_U2}) and (\ref{eq:smin_U1}) we 
deduce that 
\beq
\limsup_{n \to \infty} 
\sup_{v \in \C} \E \Big[  
\langle |\Log|, L_{\bm{V}_n^{d,v}} \rangle_0^{n^{-\delta}} 
\Big] = 0 \,. 
\label{eq:log_integrate_lower}
\eeq
Next, consider the collection of sets $\Gamma^d_{\gamma}$ 
as in (\ref{eq:gamma_d}), that corresponds to the compact 
$$
\Lambda_d := \big\{ v \in \C : |v| \in \{0,1,\ldots,d\} \big\} 
$$ 
(such that $m(\Lambda_d)=0$). In this case, $v \in \Gamma^d_{\gamma}$ 
implies that $\{ iy : y > 0 \}$ is disjoint of the set $\Gamma_\gamma^{d,|v|}$ 
of Lemma \ref{lem:induction_step}. For such values of $v$ we thus
combine the bound (\ref{eq:unif_bound}) of Lemma \ref{lem:induction_step} 
with \cite[Lemma 15]{gkz}, 
to deduce that for any integer $d \ge 1$ and $\gamma \in (0,1)$ 
there exist finite $n_0, M$ and positive $\kappa$ (depending only on 
$d$ and $\gamma$), for which 
\beq 
\E \big[ L_{\bm{V}_n^{d,v}}(-y,y) \big] \le 2 M (y \vee n^{-\kappa})
\, \qquad \forall n \ge n_0,\, y>0,\, v \in \Gamma^d_{\gamma} \,.
\label{eq:gkz-lem15}
\eeq 
Imitating the derivation of \cite[Eqn. (49)]{gkz}, 
we get from (\ref{eq:gkz-lem15}) that for some finite 
$C=C(d,\gamma,\delta)$,
any $\vep \le e^{-1}$, $n \ge n_0$ and $v \in \Gamma^d_{\gamma}$, 
\beq
\E 
\Big[ 
\langle |\Log|, L_{\bm{V}_n^{d,v}} \rangle_{n^{-\delta}}^{\vep} 
\Big] \le C \vep |\log \vep|
\,. \label{eq:log_integrate_upper}
\eeq
Thus, combining (\ref{eq:log_integrate_lower}) 
and (\ref{eq:log_integrate_upper}) we have that for any $\gamma>0$,
\beq\label{eq:vep-n-neg}
\lim_{\vep \decto 0} \limsup_{n \ra \infty} 
\sup_{v \in \Gamma^d_{\gamma}}
\E \Big[ \langle |\Log|, L_{\bm{V}_n^{d,v}} \rangle_0^{\vep} \Big] =0\,. 
\eeq
Similarly, in view of (\ref{eq:gnconv}), 
the bound (\ref{eq:unif_bound}) implies that
$$
|\Im(G_\infty^{d,v}(z))| \le M \,, 
\qquad \forall z \in \C^+ \setminus \Gamma_\gamma^{d,|v|}, \, v \in B(0,R) \,, 
$$
which in combination with \cite[Lemma 15]{gkz}, results with
$$
\Theta^{d,v} (-y,y) \le 2 M y \, 
\qquad \forall y>0,\, v \in \Gamma^d_{\gamma} 
$$ 
and consequently also
\beq\label{eq:vep-lim-neg}
\lim_{\vep \decto 0} \sup_{v \in \Gamma^d_{\gamma}}
\, \{ \langle |\Log|, \Theta^{d,v} \rangle_0^{\vep} \} = 0 \,. 
\eeq
Next, by Lemma \ref{lem:sing_val_conv}, the real valued random variables 
$X_n^{(\vep)} (\omega,v) 
:= \langle \Log  ,L_{\bm{V}_n^{d,v}} \rangle_\vep^\infty$ 
converge in probability, as $n \to \infty$, to the non-random 
$X_\infty^{(\vep)} (v) := 
\langle \Log, \Theta^{d,v}\rangle_\vep^\infty$,  
for each $v \in \C $ and $\vep>0$. This, 
together with (\ref{eq:vep-n-neg}) and (\ref{eq:vep-lim-neg}),
results with the stated convergence of (\ref{eq:log_integrates}), 
for each $v \in \Gamma^d_\gamma$, so considering $\gamma \to 0$ 
we conclude that (\ref{eq:log_integrates}) applies
for all $v \in 
\Lambda_d^c$, hence for $m$-a.e. $v$.
 
\noindent
Turning to prove (\ref{eq:function_log_integrates}), 
fix $\gamma>0$ and non-random, uniformly bounded $\phi$, 
supported within $\Gamma^d_\gamma$. Since 
$\{L_{\bm{V}_n^{d,v}}, v \in \Gamma^d_\gamma\}$ are all
supported on $B(0,\gamma^{-1} + d)$, for each fixed $\vep>0$,
the random variables 
$Y_n^{(\vep)} (\omega,v) := 
\phi(v) X_n^{(\vep)} (\omega,v) m(\Gamma^d_\gamma)$
with respect to the product law 
$\overline{\P} := \P \times m(\cdot)/m(\Gamma^d_\gamma)$ 
on $(\omega,v)$ are bounded, uniformly in $n$. Consequently, their 
convergence 
in $\P$-probability, for $m$-a.e. $v$, to $Y_\infty^{(\vep)} (v)$ 
(which we have already established), implies the 
corresponding $L_1$-convergence.  Furthermore, 
by (\ref{eq:vep-n-neg}) and Fubini's theorem,
$$
\overline{\E} [|Y^{(0)}_n - Y^{(\vep)}_n|] \le 
m(\Gamma^d_\gamma) \|\phi\|_\infty
\sup_{v \in \Gamma^d_{\gamma}}
\E [ |X_n^{(0)}(\omega,v) - X_n^{(\vep)}(\omega,v)| ] \ra 0 \,, 
$$
when $n \to \infty$ followed by $\vep \decto 0$.
Finally, by (\ref{eq:vep-lim-neg}), the non-random 
$Y_\infty^{(\vep)} (v) \to Y_\infty^{(0)}(v)$ 
as $\vep \decto 0$, uniformly over $\Gamma^d_\gamma$.
Consequently, as $n \to \infty$ followed by $\vep \decto 0$,
$$
\overline{\E} [|Y^{(0)}_n - Y^{(0)}_\infty|] \le 
\overline{\E} [|Y^{(0)}_n - Y^{(\vep)}_n|] + 
\overline{\E} [|Y^{(\vep)}_n - Y^{(\vep)}_\infty|] + 
\sup_{v \in \Gamma^d_\gamma} \{ |Y^{(0)}_\infty - Y^{(\vep)}_\infty| \} 
$$
converges to zero and in particular 
$$
\int_{\C} \phi (v) X_n^{(0)} (\omega,v) dm(v) \ra 
\int_{\C} \phi(v)  X_\infty^{(0)} (v) dm(v) \,, 
$$
in $L_1$, hence in $\P$-probability, as claimed. 
\qed

\section{Proof of Theorem \ref{thm:orthogonal}}
\label{sec:thm_ii}

\noindent
Following the proof of Proposition \ref{thm:unitary}, it suffices 
for establishing Theorem \ref{thm:orthogonal}, to extend the 
validity of Lemmas \ref{lem:sing_val_conv} and \ref{lem:log_integrate} 
in case of $S_n = \sum_{i=1}^{d'} U_n^i + \sum_{i>d'}^d O_n^i$.
To this end, recall that Lemma \ref{lem:new}(ii) applies regardless
of the value of $d'$. Hence, Lemmas \ref{lem:sing_val_conv} 
and \ref{lem:log_integrate} hold as soon as we establish 
Lemma \ref{lem:induction_start}, the bound (\ref{eq:smin_U}) 
on $s_{\min}(\bm{V}_n^{d,v})$, and the convergence 
(\ref{eq:gnconv}) for $d=1$. Examining Section 
\ref{sec:section_proofs_lemmas}, one finds that 
our proof of the latter three results applies as soon as 
$d' \ge 1$ (i.e. no need for new proofs if we start 
with $U_n^1$). 

\noindent
In view of the preceding, we set hereafter $d'=0$, namely 
consider the sum of (only) i.i.d Haar orthogonal matrices
and recall that suffices to prove our theorem when 
$d \ge 2$ (for the case of $d=1$ has already 
been established in \cite[Corollary 2.8]{MM}). Further, 
while the Haar orthogonal measure is \emph{not invariant} 
under multiplication by $e^{i \theta}$, it is not hard to 
verify that nevertheless 
$$
\lim_{n \to \infty} 
\E[\f{1}{n} \Tr \{O_n^k\}]= \E[\f{1}{n} \Tr \{(O_n^*)^k\}] = 0 \,,
$$
for any positive integer $k$. Replacing the identity 
(\ref{eq:unit-idn}) by the preceding and thereafter 
following the proof of Lemma \ref{lem:u_n_density}, we 
conclude that $\E [L_{\bm{O}_n^{1,v}}] \Ra \Theta^{1,v}$ 
as $n \to \infty$, for each fixed $v \in \C$. This yields 
of course the convergence (\ref{eq:gnconv}) of the 
corresponding Stieltjes transforms (and thereby 
extends the validity of Lemma \ref{lem:sing_val_conv} 
even for $d'=0$). Lacking the identity (\ref{eq:unit-idn}), 
for the orthogonal case we replace Lemma \ref{lem:induction_start} 
by the following. 
\begin{lem}\label{lem:induction_start_orthogonal}
The Stieltjes transform $G_n^{1,v}$ of the \abbr{ESD} 
$\E[ L_{\bm{O}_n^{1,v}}]$ is such that 
\begin{align}
   \big\{z\in \C^+: \ |\Im G^{1,v}_n (z)| \ge C \vep^{-2}\big\}
 \subset &  \Big\{E + i \eta: \eta \in (0,\vep^2), \notag\\
  E  \in \big(\pm(1 \pm |v|) - 2 \vep,& \pm(1 \pm |v|)+ 2 \vep \big) \cup  
\Big(\pm (|1 \pm v| - 2 \vep, \pm(|1 \pm v|) + 2 \vep\big) \Big\}\,, 
\notag
\end{align}
for some $C$ finite, all $\vep \in (0,1)$ and any $v \in \C$. 
\end{lem}

\noindent
{\em Proof}: 
We express $G_n^{1,v}(z)$ as the expectation of 
certain additive function of the eigenvalues of $O_n^1$, 
whereby information about the marginal distribution 
of these eigenvalues shall yield our control on
$|\Im(G_n^{1,v}(z))|$. 
To this end, set $g(z,r):=z/(z^2-r)$ for $z \in \C^+$, $r \ge 0$, 
and let $\phi(O_n^1):=\f{1}{2n} \Tr \{ (zI_{2n}-\bm{O}_n^{1,v})^{-1} \}$.
Clearly, 
\beq\label{eq:phi-idn}
\phi(O_n^1) = \f{1}{n} \sum_{k=1}^n g(z,s_k^2) \,, 
\eeq
where $\{s_k\}$ are the singular values of $O_n^1-v I_n$. 
For any matrix $A_n$ and orthogonal matrix $\wt{O}_n$, 
the singular values of $A_n$ are the same as those of 
$\wt{O}_n A_n \wt{O}_n^*$. Considering $A_n=O_n^1 - v I_n$, we
thus deduce from (\ref{eq:phi-idn}) that 
$\phi(\wt{O}_n O_n^1 \wt{O}_n^*)=\phi (O_n^1)$, 
namely that $\phi(\cdot)$ is a {\em central function} 
on the orthogonal group (see \cite[pp. 192]{agz}).

\noindent
The group of $n$-dimensional orthogonal matrices 
partitions into the classes ${\cO}^+(n)$ and ${\cO}^- (n)$ of 
orthogonal matrices having determinant $+1$ and $-1$, respectively. 
In case $n=2 \ell+1$ is odd, any $O_n \in \cO^\pm (n)$ 
has eigenvalues $\{ \pm 1,  e^{\pm i \theta_j}, j=1,\ldots, \ell\}$,
for some $\ul{\theta}=(\theta_1,\ldots,\theta_\ell) \in [-\pi,\pi]^\ell$.
Similarly, for $n=2\ell$ even, $O_n \in \cO^+(n)$ has eigenvalues
$\{e^{\pm i \theta_j}, j=1,\ldots, \ell\}$, whereas 
$O_n \in \cO^-(n)$ has eigenvalues
$\{-1,1,e^{\pm i \theta_j}, j=1,\ldots, \ell-1\}$. Weyl's formula 
expresses the expected value of a central function of 
Haar distributed orthogonal matrix in terms of the joint 
distribution of $\underline{\theta}$ under the probability
measures $\P^{\pm}_n$ corresponding to the classes 
$\cO^+(n)$ and $\cO^-(n)$. Specifically, it yields the expression 
\begin{align}
G_n^{1,v}(z) = \E[\phi(O_n^1)] &= 
\f{1}{2} \E^+_n[\phi(\diag(+1,R_\ell(\ul{\theta}))]+ 
\f{1}{2} \E^-_n[\phi(\diag(-1,R_\ell(\ul{\theta}))]\,, 
\quad \mbox { for }  \; n=2\ell+1, 
\notag 
\\
&=
\f{1}{2} \E^+_n[\phi(\diag(R_\ell(\ul{\theta}))]+ 
\f{1}{2} \E^-_n[\phi(\diag(-1,1,R_{\ell-1} (\ul{\theta}))]\,, 
\quad \mbox { for }  \; n=2\ell,
\label{eq:n2l}
\end{align}
where $R_\ell(\ul{\theta}):= 
\diag(R(\theta_1), R(\theta_2), \cdots, R(\theta_\ell))$ 
for the two dimensional rotation matrix
\beq
R(\theta)=\begin{bmatrix}
\cos\theta & \sin\theta \\
-\sin\theta & \cos\theta
\end{bmatrix}
\notag
\eeq
(see \cite[Proposition 4.1.6]{agz}, which also provides
the joint densities of $\ul{\theta}$ under $\P^\pm_n$).

In view of (\ref{eq:phi-idn}) and 
(\ref{eq:n2l}),
to evaluate $G_n^{1,v}(z)$ we need the singular values of 
$R_\ell(\ul{\theta})-v I_\ell$. Since this 
is a block-diagonal matrix, its singular values are those
of the $2 \times 2$ block diagonal parts 
$R(\theta_j)-vI_2$ for $1 \le j \le \ell$.
Setting $v:=|v| e^{i \psi}$
it is easy to check that the singular values of $R(\theta) - vI_2$ 
are precisely square-root of the eigenvalues of  
$(1+|v|^2)I_2 - |v| (e^{-i \psi} R(\theta) + e^{i \psi} R^*(\theta))$,
which turn out to be $1+|v|^2 - 2|v| \cos(\theta \pm\psi)$. 
Combining this with (\ref{eq:phi-idn}) and (\ref{eq:n2l}) we 
obtain in case $n=2\ell+1$, that
\begin{align}
G_n^{1,v}(z) = \f{1}{2n} \Big\{ & g(z,|1 -v|^2) +
\sum_{k=0}^1 \sum_{j=1}^\ell 
\E^+_n [g(z,1+|v|^2 - 2|v| \cos(\theta_j 
+ (-1)^k \psi))] \notag \\
+ & g(z,|1 +v|^2) +
\sum_{k=0}^1 \sum_{j=1}^\ell \E^-_n [g(z,1+|v|^2 - 2|v| 
\cos(\theta_j + (-1)^k \psi))] \Big\} \,.
\label{eq:simplify_ortho_1}
\end{align}
The same expression applies for $n=2\ell$, except for having
the latter sum only up to $j=\ell-1$. 
Next, recall that under $\P^\pm_n$ the random variables 
$\{\theta_j\}$ are exchangeable, each having the same
density $q^{\pm}_n(\cdot)$ which is bounded, uniformly in $n$
(see the diagonal terms in \cite[Proposition 5.5.3]{forrester}; 
for example, $q^{\pm}_{2\ell+1} (\theta) = 
\frac{1}{2 \pi}( 1 \mp \sin(2 \ell \theta)/(2 \ell \sin \theta))$,
is bounded by $1/\pi$, uniformly over $\theta$ and $\ell$). 
Further, $g(z,r) \in \C^-$
for all $r \ge 0$ and $z \in \C^+$. Hence, 
for some $C$ finite, all $n \ge 3$, $v \in \C$ and $z \in \C^+$,
\begin{align}
|\Im(G_n^{1,v}(z))| \le& \f{1}{2n} |\Im (g(z,|1 -v|^2))| 
+ \f{1}{2n} |\Im (g(z,|1 +v|^2))| \notag \\
&+ C \Big| 
\Im\big\{ \f{1}{2\pi} \int_{-\pi}^{\pi} g(z,1+|v|^2 - 2|v| \cos(\theta \pm\psi)) 
d \theta \big\} \Big| \,.
\label{eq:ortho_bd}
\end{align}
The last expression in (\ref{eq:ortho_bd}) does not depend on 
$\pm \psi$ and is precisely the imaginary part of the Stieltjes transform 
of the symmetrization of the probability measure $|e^{i \theta} - |v||$, 
where $\theta \sim \dU(0,2 \pi)$. While proving Lemma \ref{lem:u_n_density} 
we saw that the expected \abbr{ESD} of $\bm{U}_n^{1,v}$ has 
the latter law, 
hence 
the conclusion of Lemma \ref{lem:induction_start} 
applies for the last expression in (\ref{eq:ortho_bd}).
To complete the proof, simply note that 
$\Im (g(E+i\eta,s^2)) \le 1$ as soon as $|E \pm s| \ge \sqrt{\eta}$
(and consider $s = |v \pm 1|$). 
\qed

\medskip
Now, using Lemma \ref{lem:induction_start_orthogonal} 
for the basis $d=1$ of an induction argument 
(instead of Lemma \ref{lem:induction_start}), 
and with Lemma \ref{lem:new}(ii) serving
again for its inductive step, we obtain here
the same conclusion as in 
Lemma \ref{lem:induction_step}, except for 
replacing $\Gamma_\gamma^{d,|v|}$ by the 
union $\wt{\Gamma}_\gamma^{d,v}$ of open balls
of radius $\gamma$ centered at the points 
$\pm m \pm 1 \pm |v|$ \emph{and} 
$\pm m \pm |1 \pm v|$ for $m=0,\ldots,d-1$. 
Turning to prove Lemma \ref{lem:log_integrate}, 
this translates to taking in this case the
sets $\Gamma^d_\gamma$ which correspond 
via (\ref{eq:gamma_d}) to the compact 
$$
\Lambda_d := \big\{ v \in \C : |v| \in \{0,1,\ldots,d\}, 
\quad \mbox { or } \quad 
|v \pm 1| \in \{0,1,\ldots,d-1\}\, \big\} 
$$ 
(of zero Lebesgue measure), thereby assuring 
that $\{i y: y > 0 \}$ is disjoint of 
$\wt{\Gamma}_\gamma^{d,v}$ whenever 
$v \in \Gamma^d_\gamma$. One may then easily 
check that the proof of Lemma \ref{lem:log_integrate} 
(and hence of the theorem),
is completed upon establishing the following 
weaker form of (\ref{eq:smin_U}).
\begin{lem}\label{lem:log_integrate_lower_ortho}
For some $c_1>0$, $c_2 <\infty$, the sum $S_n$ of $d \ge 2$ 
independent Haar orthogonal matrices and any $\gamma \in (0,1)$, 
there exist $C'=C'(d,\gamma)$ finite and events $\{\cG_n\}$ determined
by the minimal and maximal singular values of $S_n$, such that 
$\P(\cG_n^c) \ra 0$ as $n \to \infty$, and for any $n, t \ge 0$,
\beq
\sup_{v \in \Gamma_\gamma^d}
\P\Big(\cG_n \cap \{s_{\min}(\bm{V}_n^{d,v}) \le t\}\Big) 
\le C' t^{c_1} n^{c_2} \,. 
\label{eq:smin_O}
\eeq
\end{lem}
\noindent 
{\em Proof}: We use here \cite[Theorem 1.3]{vershynin} (instead of 
\cite[Theorem 1.1]{vershynin} which applies only for Haar
unitary matrices), 
and introduce events $\cG_n$ under which the condition 
\cite[Eqn. (1.2)]{vershynin} holds. Specifically, let 
$D_n = \diag(r_1,r_2,\ldots,r_n)$ denote the diagonal matrix of 
singular values of $S_n$, ordered so that $r_1 \ge r_2 \ge \ldots \ge r_n$ 
and 
\beq
\cG_n := \{r_n \le \f{1}{2} \quad \mbox { and } \quad r_1 \ge 1 \} \,.
\notag
\eeq
Let $O_n$ be Haar distributed $n$-dimensional orthogonal matrix, 
independent of $\{O_n^i, i=1,\ldots,d\}$, noting that
$O_n$ is independent of $-O_n S_n$, with the latter having 
the same law and singular values as $S_n$. 
Further, the singular values of $\bm{V}_n^{d,v}$ equal to 
those of $vI_n - S_n = O_n^*(v O_n - O_n S_n)$, hence for 
any $n$ and $t \ge 0$,
\begin{align}
q_{n,v}(t) := \P\Big(\cG_n \cap \{s_{\min}(\bm{V}_n^{d,v}) \le t\}\Big) 
&= \P\Big(\cG_n \cap \{s_{\min}(v O_n + S_n) \le t\} \Big)\notag.
\end{align}
Next, by the singular value decomposition 
$S_n=(O_n^{'})^* D_n (O_n^{''})^*$ for some 
pair of orthogonal matrices $O_n'$ and $O_n^{''}$. 
Conditional on $D_n$, $O_n^{'}$ and $O_n^{''}$, 
the matrix $O_n^{'} O_n O_n^{''}$ is again Haar 
distributed, hence independent of $D_n$ (and of $\cG_n$). 
Consequently, for any $v \ne 0$,  
\begin{align}
q_{n,v} (t) 
= \P\Big(\cG_n \cap \{s_{\min}(v O_n^{'} O_n O_n^{''} + D_n) \le t\}\Big)
= \P\Big(\cG_n \cap \{|v| s_{\min}(O_n + v^{-1} D_n) \le t\}\Big). \notag
\end{align}
Now from \cite[Theorem 1.3]{vershynin} we know that for some 
absolute constants $c_1>0$ and $c_2<\infty$, 
\beq
\P(|v| s_{\min} (O_n + v^{-1} D_n) \le t \,|\, D_n) 
\le \Big(\f{t}{|v|}\Big)^{c_1} 
\Big(\f{K n}{\delta}\Big)^{c_2}\,, 
\label{eq:rv_orthogonal}
\eeq
provided \cite[Eqn. (1.2)]{vershynin} holds for $v^{-1} D_n$, 
some $K \ge 1$ and $\delta \in (0,1)$. That is, when
\beq
r_1 \le K |v| \text{ and } r_1^2 \ge r_n^2 + \delta |v|^2 \,.
\label{eq:rv_orthogonal_check}
\eeq
In our setting the singular values of $S_n$ are uniformly bounded by
$d$ and $|v| \in (\gamma,\gamma^{-1})$ throughout 
$\Gamma_\gamma^d$. Hence, the event $\cG_n$ implies that 
(\ref{eq:rv_orthogonal_check}) holds for $K=d/\gamma$ 
and $\delta=\gamma^2/2$. Thus, multiplying both sides of 
(\ref{eq:rv_orthogonal}) by $\bI_{\cG_n}$ and 
taking the expectation over $D_n$ yields the inequality
(\ref{eq:smin_O}) for some finite $C'=C'(d,\gamma)$.

\noindent
Proceeding to verify that $\P(\cG_n^c) \to 0$ as $n \to \infty$, 
recall \cite[Proposition 3.5]{haagerup_larsen} that $\Theta^{d,0}$
is the symmetrization of the law $\mu_{|s_d|}$, for the sum 
$s_d=u_1+\cdots+u_d$ of $\star$-free Haar unitary 
operators $u_1,\ldots,u_d$, and 
\cite[Eqn. (5.7)]{haagerup_larsen} that for $d \ge 2$ 
the measure $\mu_{|s_d|}$ on $\R^+$ has the density 
\beq
\f{d\mu_{|s_d|}}{dx} = \f{d\sqrt{4(d-1)-x^2}}{\pi(d^2-x^2)} \bI_{[0,2\sqrt{d-1}]}(x) \,,
 \label{eq:limit_sv}
\eeq
so in particular both $\mu_{|s_d|}((0,1/2))$ and 
$\mu_{|s_d|}((1,3/2))$ are strictly positive. 
Further, from Lemma \ref{lem:sing_val_conv}
we already know that the symmetrization of the \abbr{ESD} 
$\nu_{|S_n|}$ of $D_n$, converges weakly, in probability, 
to $\Theta^{d,0}$ and consequently, $\nu_{|S_n|}$ 
converges weakly to $\mu_{|s_d|}$, in probability. From 
the preceding we deduce the existence of $g \in C_b(\R^+)$ 
supported on $[0,1/2]$, such that $\langle g, \mu_{|s_d|} 
\rangle \ge 1$ and that for such $g$,  
\beq
\P(r_n > 1/2) \le \P( \langle g, \nu_{|S_n|} \rangle =0) 
\le \P \Big(|\langle g, \nu_{|S_n|} \rangle  
- \langle g, \mu_{|s_d|} \rangle | >1/2 \Big) \ra 0\,, \label{eq:s_min}
\eeq
as $n \to \infty$. Similarly, considering $g \in C_b(\R^+)$ 
supported on $[1,3/2]$ for which $\langle g, \mu_{|s_d|} \rangle \ge 1$,
we get that $\P(r_1  < 1) \ra 0$, from which we conclude that 
$\P(\cG_n^c) \to 0$.  
\qed

\section{Proof of Proposition \ref{rmk:gkz_improved}}
\label{sec:prop}

\noindent

The main task here is to show that for $m$-a.e.$v \in \C$,
the logarithm is uniformly integrable with respect to the 
\abbr{ESD} of $|U_n T_n - v I_n|$. As shown in \cite{gkz}, 
setting $\rho = |v|$, this is equivalent to such uniform 
integrability for the \abbr{ESD} $\nu_n^v$ of the matrix
$\bm{Y}_n^v$ (per (\ref{eq:Y_n_T_n})).
The key for the latter is to show that $\Im(G_n(\cdot))$ 
is uniformly bounded on $\{i \eta: \eta > n^{-\kappa_1} \}$ for 
some $\kappa_1 >0$ and Lebesgue almost every $\rho$ 
(see proof of \cite[Proposition 14 (i)]{gkz}). In \cite{gkz}, this was done 
under the assumption of \cite[Eqn. (3)]{gkz}, whereas here we 
show that the same holds under the weaker condition (\ref{eq:modi_T_n}).  

\noindent
To this end, \cite[Lemma 10]{gkz} yields (analogously to 
Lemma \ref{lem:sing_val_conv}), the weak convergence, 
in probability, of $\nu_n^v$ to $\nu^v$, as well as the 
identities and bounds \cite[Eqn. (34)--(38)]{gkz}, 
without ever using \cite[Eqn. (2) or Eqn. (3)]{gkz}. 
The same applies for \cite[Lemma 11 and Lemma 12]{gkz} 
which validate the Schwinger-Dyson equation 
\cite[Eqn. (38)]{gkz} for all $n$ large enough,
any $\Im(z) > C_1 n^{-1/4}$ and $\rho \in (0,R]$. 
We then use Lemma \ref{lem:new}(i) to bypass \cite[Lemma 13]{gkz}. 
Specifically, from (\ref{eq:modi_T_n}), using 
Lemma \ref{lem:new}(i) we have that 
for every $\vep \in (0,1/2)$ and finite $R$, there 
exist finite $M_1$ and $\kappa_1>0$ 
depending only on $R$ and $\vep$ such that 
for every $\rho \in [R^{-1},R]$,
\beq\label{eq:bd-styn}
\{z: \Im (z) > n^{-\kappa_1}, |\Im(G_{n}(z))| > M_1 \} 
\subset \Gamma_{2\vep}^\rho
\eeq
where $\Gamma_\gamma^\rho$ denotes the union of open 
balls of radius $\gamma>0$ centered at points from 
the symmetric subset $K \pm \rho$ of $\R$. 
Having (\ref{eq:bd-styn}) instead of the bound 
(\ref{eq:unif_bound}) of Lemma \ref{lem:induction_step},
we consider here the closed set 
$\Lambda_K := \{ v \in \C : |v| \in K \}$ such that
$m(\Lambda_K)=0$, the bounded, open sets $\Gamma_\gamma$, $\gamma>0$,
associated with $\Lambda_K$ via (\ref{eq:gamma_d}), and the 
corresponding collection $\cS \subset C_c^\infty(\C)$ of 
test functions. Using this framework and following 
the proof of Lemma \ref{lem:log_integrate}, we deduce 
that
\beq
\langle \Log, \nu_n^v \rangle \ra \langle \Log, \nu^v \rangle, \notag
\eeq
in probability for each $v \in \Gamma_\gamma$, and consequently
for $m$-a.e. $v \in \C$. Then, utilizing our assumption 
(\ref{eq:assumption_1}) on the uniformly bounded support 
of the relevant \abbr{ESD}-s, we have that further,
for any fixed $\phi \in \cS$, 
\beq
\int_\C \phi(v) \langle \Log, \nu_n^v \rangle d m(v) \ra 
\int_\C \phi(v) \langle \Log, \nu^v \rangle d m (v)\,,
\notag
\eeq
in probability. Since $\Theta$ not a Dirac measure, we know 
from \cite[Theorem 4.4]{haagerup_larsen} and \cite[Remark 8]{gkz} 
that $\mu_A$ has a density with respect to the Lebesgue measure 
on $\C$. Consequently $\mu_A(\Lambda_K)=0$, and following the
same argument as in the proof of Proposition \ref{thm:unitary},
we get part (a) of Proposition \ref{rmk:gkz_improved}.

\noindent 
For parts (b) and (c) of the proposition see \cite[Remark 8]{gkz} 
(which does not involve \cite[Eqn. (2) or Eqn. (3)]{gkz}). For
part (d) recall that Lemma \ref{lem:new}(i) applies even 
in case $U_n$ is replaced by a Haar distributed orthogonal 
matrix $O_n$, as does the relevant analysis 
from \cite{gkz} (c.f. proof of \cite[Theorem 18]{gkz}). 
Hence, following the same argument as in the unitary case, 
the proof is complete once we establish the analog of 
Lemma \ref{lem:log_integrate_lower_ortho}.
That is, specify events $\cG_n$ determined by $T_n$, 
such that $\P(\cG_n^c) \to 0$ as $n \to \infty$ and 
\beq\label{eq:bd-smin6}
\sup_{v \in \Gamma_\gamma} 
\P\Big( \cG_n \cap \big\{ 
|v| s_{\min} (O_n + v^{-1} T_n) \le t \big\} \Big) 
\le C' t^{c_1} n^{c_2} \,,
\eeq
for any $\gamma>0$, some $C'=C'(\gamma)$ finite and all $t$, $n$. 
To this end, with $\Theta$ non-degenerate, there exist 
$\xi>0$ and $b_+^2 \ge b_-^2 + \xi$, such that both 
$\Theta([0,b_-))$ and $\Theta((b_+,M])$ are positive. 
Consequently, setting $T_n = \diag(r_1,\ldots,r_n)$ 
with $r_1 \ge r_2 \ge \ldots \ge r_n$, it follows from the 
weak convergence of $L_{T_n}$ to $\Theta$ (in probability), 
that $\P(\cG_n^c) \to 0$ 
for $\cG_n := \{ r_n \le b_-$ and $r_1 \in [b_+,M] \}$ (by the
same reasoning as in the derivation of (\ref{eq:s_min})). 
Further, (\ref{eq:bd-smin6}) follows by an application of 
\cite[Theorem 1.3]{vershynin} conditional upon $T_n$
(where \cite[Eqn. (1.2)]{vershynin} holds under $\cG_n$ 
for $v^{-1} T_n$, $v \in \Gamma_\gamma$, $K=M/\gamma$ and 
$\delta = \xi \gamma^2$, see  
(\ref{eq:rv_orthogonal})-(\ref{eq:rv_orthogonal_check})).

\vskip 5pt

\noindent
{\bf Acknowledgment} 
We thank Ofer Zeitouni for suggesting to look at the sum of $d$ 
unitary matrices and for pointing our attention to preliminary 
versions of \cite{bordenave, vershynin}. 
We also thank the anonymous referees whose valuable comments 
on an earlier draft helped improving the presentation of this 
work.

\end{document}